\documentclass[a4paper,12pt]{article}
\usepackage{mathptmx}
\usepackage{amsthm,amsfonts, indentfirst, latexsym, bm,graphicx,amsmath,amssymb,mathrsfs,verbatim,hyperref,pdfpages}
\usepackage{caption}
\usepackage{subcaption}
\usepackage{color}
\usepackage[T1]{fontenc}
\usepackage{dsfont}
\usepackage{latexsym}
\usepackage{url}
\usepackage{array}
\usepackage{tikz}
\usepackage{pgfplots}
\usepackage{graphicx}
\usepackage{layout}
\usepackage[titletoc,title]{appendix}% this is for appendix
\usepackage{bbm} % this is for indicator function
\usepackage{algorithm,algorithmic}
\numberwithin{equation}{section}
\newtheorem{theorem}{Theorem}[section]
\newtheorem{corollary}[theorem]{Corollary}
\newtheorem{proposition}[theorem]{Proposition}

\newtheorem{lemma}[theorem]{Lemma}

\usepackage{chngcntr}
\counterwithout{equation}{section} % undo numbering system provided by phstyle.cls
\usetikzlibrary{arrows, automata,positioning,calc,shapes,decorations.pathreplacing}
  \colorlet{greencolor}{green!50!black}
  \colorlet{textcolor}{red}
  \colorlet{tancolor}{orange!80!black}
  \colorlet{bluecolor}{blue}
\usepackage{amsmath,amsfonts,amssymb,amsthm,epsfig,epstopdf,url,array}

\definecolor{plotcolor1}{rgb}{0,0.447,0.741}
\definecolor{plotcolor2}{rgb}{0.741,0,0.447}
\definecolor{plotcolor3}{rgb}{0,0.741,0.294}
\definecolor{plotcolor4}{rgb}{0.741,0.294,0}
\definecolor{plotcoloraux}{rgb}{0.447,0.447,0.447}
\tikzset{plotstyle1/.style={color=plotcolor1,solid,line width=1.0pt}}
\tikzset{plotstyle2/.style={color=plotcolor2,densely dashed,line width=1.0pt}}
\tikzset{plotstyle3/.style={color=plotcolor3,dotted,line width=1.0pt}}
\tikzset{plotstyle4/.style={color=plotcolor4,loosely dashed,line width=1.0pt}}
\tikzset{auxlines/.style={color=plotcoloraux,solid,line width=0.5pt}}
\newcommand{\markersize}{0.7pt}
\tikzset{discretemarkers/.style={mark=*,mark options={solid},mark size=\markersize}}
\usepackage{makecell} %for changing the lines in table
\def\VR{\kern-\arraycolsep\strut\vrule &\kern-\arraycolsep}%for put sth on the above or left side of matrices
\def\vr{\kern-\arraycolsep & \kern-\arraycolsep}
\title{ A Finsler Geometrical Programming  for the Nonlinear Complementarity Problem of Traffic Equilibrium}

\author{Azam Asanjarani\footnote{The University of Auckland}}

\date{ }
\begin{document}
\maketitle

%%%%%%%%%%%%%%%%%%%%%%%%%%%%%%%%%%%%%%%%%%%%%%%%
%%%%%%%%%%%%%%%%%%%%%%%%%%%%%%%%%%%%%%%%%%%%%%%%
\begin{abstract}
%Geometry as one of the oldest sciences contributes greatly to different branches of mathematics including optimization.

This work is a geometrical approach to the optimization problem  motivated by transportation system management. First, an attempt has been made to furnish  a comprehensive account of geometric programming based on the elementary Finsler geometry in $I\!\!R^n$. Then, a Finslerian dynamical model for the nonlinear complementarity problem of traffic
equilibrium  is presented that
 may be applied for  different equilibrium problems.
\end{abstract}

%\begin{keywords}Nonlinear complementarity,Traffic 
%Equilibrium, Randers, Finsler
%\end{keywords}
%\begin{AMS} 90C33, 58B20
%\end{AMS}
%\pagestyle{myheadings}
%\thispagestyle{plain}
%\markboth{PAUL DUGGAN}{USING SIAM'S \LaTeX\ MACROS}

\section{Introduction}

 Finsler metrics, as a natural  extension of Riemannian and Euclidean metrics, becomes more present in mathematical programming. Particularly,
Finsler metrics originated in the study of calculus of variation by P.
Finsler in the first quarter of twentieth century, is more reliable
in optimization theory than other two metrics, see for instance \cite{antonelli2005differential, antonelli2013theory,  antonelli2005finslerian, kielanowski2019functional, kristaly2008optimal}.
In this paper, we restrict
our attention to the Finsler metrics defined on  $I\!\!R^n$ in order to formulate the  problems of nonlinear complementarity problem  (NCP) and traffic equilibrium with geometric  models.

In a general sense, an optimisation problem or a mathematical programming mainly refers to choosing the best (optimal) elements from a set of available alternatives. If the set includes different paths between two fixed points, we may be able to  consider a metric on the set and define the optimised paths as the shortest paths or geodesics of this metric. Definition and  axioms of the geometry of a space  basically depend on the definition of the shortest path in that space.
 For a well-known metric (e.g. Euclidean, Riemannian, or Finslerian metric),  this happens by establishing the relevant geometry to have a better understating of its geometric behaviour. This phenomenon may be called the geometrisation of mathematical programming. Here, we present some examples of  geometrisation of  complementarity problems and traffic equilibriums. This work may serve as the first stage of applying  Finsler metrics in these types of mathematical programming.
 
The need for measures to reduce congestion and manage traffic demand  in the metropolitan
areas is becoming more serious as the population grows. A
functioning society depends on the mobility provided by the
transportation network to enable its members to participate in
essential activities such as production, consumption, communication,
and recreation. However, it is necessary for society also  to
introduce congestion-relief measures for the quality of life, the
environment, and the safety of the citizens not to deteriorate.

 Any well-founded traffic model recognises the individual network user's
right to decide about when, where, and how to travel  \cite{marcotte2007traffic}. However, there is always a
 conflict of interest in the transportation system: A typical traveller 
expected to choose an optimal route, the route that minimises 
 the combined travel time and cost,  given the network conditions during the travel time (the Wardrop's  principle). But, society's goal often is to reduce the average travel times and the environment's damage. The aggregate effect of these
decisions is a network flow which does not minimise the total system
costs. The Wardrop conditions are
frequently used in urban traffic planning to predict and analyse
traffic flows. The equilibrium conditions are based on the
assumption of rational route choice behaviour by  individual users
and define a situation where the travellers' routes are
 the shortest paths between the given set of origins and destinations.

 In the design
and management of urban transportation systems, there is a need for
efficient planning tools for analysing and predicting future scenarios. Such tools rely on
the mathematical model of the transportation system, the
users of the system and their aggregate behaviour under different
traffic conditions \cite{bellomo2002modelling}.
 From the mathematical point  of view, it is assumed that the area under study is strongly connected and
there is at least one path between any origin-destination pair. Furthermore, 
the transportation system is represented by a network that
describes  transportation possibilities (i.e., roads, transit
lines, etc.) between origin-destination pairs.

The nonlinear complementarity problem (NCP)
 is an important problem in mathematical programming that has many applications in different fields, see
for instance \cite{harker1990finite}. One of the well-known NCP functions is the 
Fischer-Burmister function \cite{fischer1997solution} that is used extensively in the solution of nonlinear
complementarity and variational inequality problems \cite{facchinei2007finite}. In
1999, Qi in \cite{qi1999regular} showed that the  and
its several variants, such as the Tseng-Luo NCP function (\cite{luo1997new})
and the Kanzow-Kleinmichel NCP function (\cite{kanzow1998new}) are smooth in the areas away from the origin and are strongly semi-smooth at the origin. The
Fischer-Burmister function and its variants are also irrational. Moreover, Qi proposed  a  class of piecewise rational NCP functions
with the same strongly semi-smooth property, see \cite{qi1999regular}.

Here, we study the traffic equilibrium problem from the Finsler geometrical
point of view and show that the shortest paths of traffic problem
are solutions of a system of nonlinear equations known as geodesics of a Randers metric. Also, we show that for a dynamical transition network, optimised routes  are  minimal of an integration of the Fischer-Burmister function. This leads to a Finsler geometrical model for the  dynamical Wardrop's user equilibrium problem.

The remainder of this paper is structured as follows. In Section~\ref{sec:1}, we briefly define a Finsler metric on $I\!\!R^n$, the NCP and the traffic equilibrium problem. In Section~\ref{Sec:Finsler-Traffic-model} we present a
Finsler geometrical model  for the traffic problem and show that the optimal routs (shortest paths) of traffic problem can be presented as 
geodesics of a Randers metric.
In Section~\ref{Sec:Ward}, first we indicate that the Fischer-Burmister function can be considered as a special Randers metric.  Then, we
present a dynamical model for solving NCP. Finally, we propose a mathematical model for solving the dynamical
Wardrop's user equilibrium problem. These
results can apply for different equilibrium problems.

%%%%%%%%%%%%%%
%%%%%%%%%%%%%%
%%%%%%%%%%%%%%%

\section{Preliminaries}\label{sec:1}
In this section,  we briefly introduce  Finsler metrics on $I\!\!R^n$, the nonlinear complementarity problem (NCP) and the traffic equilibrium problem. We assume that the reader has no background knowledge of this topic.
\subsection{Finsler metrics on $I\!\!R^n$}

Let $M=I\!\!R^n$ be the real n-dimensional space.  The set of all tangent vectors at the point $x\in M$ is called the tangent space of $M$ and denoted by $T_xM$. The set of all tangent spaces at $x\in M$ is called the tangent bundle of $M$ and denoted by $TM$, i.e. $\displaystyle     TM:=\bigcup_{x\in M} T_xM$.  Each element of $TM$ has the form $(x,y)$ consisting of the point $x\in M$ and the tangent vector $y \in T_xM$ at the point $x$. Let   the natural projection  $\pi:TM\rightarrow M$  be defined as   $\pi(x,y):=x$. The pair $(M,\overline{F})$  is said to be a \emph{Finsler space} where $\overline{F}:TM \rightarrow
[0,\infty )$   is a function with the following properties:
\begin{enumerate}
\item $\overline{F}(x,y)$ is
differentiable on $TM \backslash \{0\} $.
\item $\overline{F}$ is positively
homogeneous of degree one in $y$, i.e. $ \overline{F}(x,\lambda
y)=\lambda \overline{F}(x,y),  \forall\lambda>0$.
\item The Hessian matrix of
$\overline{F}^{2}$, 
\begin{equation}
\label{eq:g}
g=(g_{ij}):=\left({1 \over 2} \left[
\frac{\partial^{2}}{\partial y^{i}\partial y^{j}} \overline{F}^2 
\right]\right),
 \end{equation}
 is positive-definite on $TM \backslash \{0\}  $ where  $i, j \in 1, \cdots, n$.
\end{enumerate}
The function $\overline{F}$   is called a \emph{Finsler  structure} and $g$ is called its associated  \emph{Finsler metric} on manifold $M$. Here, we  use the notation  $\overline{F}$ for the Finsler structure, rather than the usual notation $F$, to avoid any confusion with the function $F$ in the nonlinear complementarity problem.

  We can define a norm function on the vector space $T_xM$ based on  any Finsler metric $g$ on $M$.
 Let $\{U,(x^i)\}$ or simply $(x^1,...,x^n)=(x^i): U \rightarrow I\!\!R^n$ be a local coordinate system on an open subset $U\subset M$ around the point $x \in M$. The coordinate system $\{U,(x^i)\}$ induces a coordinate basis $\frac{\partial} {\partial x^i}$ on $T_xM$. Hence, we can write each tangent vector in the form $y=y^i \frac{\partial} {\partial x^i}\,$, where we apply  the \emph{Einstein summation convention or Einstein notation}, that is, whenever an index variable appears twice in a single term, once in an upper (superscript) and once in a lower (subscript) position, it implies that we are summing over all of its possible values.

 Let $x$  be a point of $M$ with the local coordinate system $(x^i)$, then it generates  a local coordinate system $(x^i, y^i)$ on $TM$. The pair  $(x^i, y^i)$ consisting of the position  $x\in M$ and the direction  $y\in T_xM$, is known as the line element of $(M,\overline{F})$.
The Finsler structure $\overline{F}(x,y)$ is said to be \emph{Riemannian} if it is independent of the  direction $y$ and by homogeneity of $\overline{F}$ it is equivalent to say that   $\displaystyle\frac{\partial \overline{F}}{\partial y^i}=0$ for $i \in 1, \cdots, n$. The Finsler structure $\overline{F}(x,y)$ is said to be  \emph{Euclidean} if $g_{ij}(x,y)=\delta_{ij}$,  where $g_{ij}$ is defined  by  \eqref{eq:g} and $\delta_{ij}$ is the Kronecker delta.

%%%%%%%%%%%%%%%%%%%
%%%%%%%%%%%%%%%%%%%
%%%%%%%%%%%%%%%%%%%
We denote by $\Gamma(p,q)$ the collection of all piecewise
$C^\infty$ curves $\sigma:[a,b]\subset I\!\!R \longrightarrow M$ on
$(M,\overline{F})$ with $\sigma(a)=p$ and $\sigma(b)=q$. Then the
length of $\sigma$ is defined by
$$
J(\sigma):=\int_a^b\overline{F}(\sigma,\frac{d\sigma}{dt})dt,
$$
where
$\frac{d\sigma}{dt}=\frac{d\sigma^i}{dt}\frac{\partial}{\partial
x^i}\in T_{\sigma(t)}M$. Define the map
$$\begin{array}{clc}
d:M\times M\longrightarrow [0,\infty)\\
d(p,q):=\inf_{\Gamma(p,q)} J(\sigma).
\end{array}
$$
It can be shown that $(M,d)$ satisfies the axioms of a metric space,
expect the symmetry property, see  \cite{antonelli2013theory} and \cite{shen2001lectures}.
In the original sense, a geodesic is a generalization of the notion of a "straight line" in a Euclidean space. However, when the metric is not Euclidean, geodesics are not necessarily straight lines and define as follows. A piecewise $C^\infty$ curve $\sigma:[a,b]\longrightarrow M$ with
$\sigma(a)=p$ and $\sigma(b)=q$ on the space  $M$ with the Finsler structure $\overline{F}$
is said to be a \emph{geodesic } if it is a  minimal curve in $\Gamma(p,q)$,
 with a constant velocity, i.e. $J(\sigma)=d(p,q)$. Hence,  geodesics are defined to be (locally) the shortest path between points in the space.

For any Finsler metric $g$ on  $M$, we can define  an inner product $g(.,.):=<.,.>_g$ on the tangent space $T_xM$ (here, for $M=I\!\!R^n$ we have ,  $T_xM\cong I\!\!R^n$). 
In the local coordinate system  $\{U,(x^i)\}$, for all $ y,z \in T_xM$ we have $g(y,z)=g_{ij}y^iz^j$.
Each inner product defines a norm $|| y||_g:=<y,y>_g=g_{ij} y^i y^j$
for a vector $y\in T_xM$ with respect to $g$. Hence,  a vector $y$ on the tangent space $T_xM$ can have different  norms according to different Finsler metrics defined on $M$. 
For a Riemannian metric  $g_{ij}(\cdot)$, 
% function
%$\overline{F}:U\times I\!\!R^3\longrightarrow I\!\!R$, defined by
$\overline{F}(x^i,y^i)=(g_{ij}y^iy^j)^{\frac{1}{2}} + b_i(x) y^i$,
where $\beta= b_i(x)y^i$ is a  1-form on $M $ with $|| \beta||_g <1$, is a Finsler
structure on $M$ and its associated Finsler metric is called a  \emph{Randers metric}, see \cite{shen2001lectures} for more details.
%%%%%%%%%%%%%%%
%%%%%%%%%%%%%%%%%
%%%%%%%%%%%%%%%%%
%{\color{red}
% If the norm $|| . ||_g$ on the tangent space $T_xM$ is  related to a mathematical programming, then the metric $g$ on $M$ determines the shortest paths as geodesics of $g$.}
%\subsection{ Geometrization of some Mathematical Programs}
%%\paragraph{}

  \subsection{Nonlinear Complementarity Problem}
 The \textit{nonlinear complementarity problems} (NCPs)  arise from optimisation theory, engineering and
economic applications. The restricted NCP functions which can be
used to reformulate NCPs as
constrained optimisation problems are first introduced by Yamashita
\cite{yamashita1998properties}. Furthermore, the discretisation of infinite-dimensional
variational inequalities leads to linear and nonlinear
complementarity problems with many degrees of freedom. Luo and Tseng
\cite{luo1997new} proposed a class of merit functions for the NCPs and showed that they satisfy several
interesting properties under some assumptions. Kanzow,
Yamashita and Fukushima have used the similar idea and proposed  new merit functions for the NCP in \cite{kanzow1997new}.

 For a given smooth mapping $F:I\!\!R^n
\rightarrow I\!\!R^n$, the NCP consists of
finding a vector $x\in I\!\!R^n$ such that
\begin{equation}\label{3}
x\geq 0, \qquad F(x)\geq 0, \qquad F(x)^T.x=0,
\end{equation}
 where $F(x)^T$ is the transpose of $F(x)$.
 Complementarity problems can be reformulated as systems of
nonlinear equations in several ways. A large number of methods has been developed based
on \textit{Newton method} and its generalisations. In order to overcome some disadvantages of the
class of non-smooth Newton methods and the class of the so-called
smoothing methods, a third class known as \textit{Jacobian smoothing
methods}  developed in \cite{kanzow1996equivalence}.  A Jacobian smoothing method is a mixture of Newton and non-smooth methods and is  based  on solving  the mixed Newton equation: a  combination of the original semi-smooth
function and the Jacobian of the smooth operator of the \textit{Clarke generalised Jacobian}  \cite{clarke1990optimization}.
 There are also several  related algorithms entitled
\textit{Newtonian algorithms}, see  for instance  \cite{ facchinei1997new, jiang1999global, kanzow1996equivalence, qi2000smoothing, yamashita1998properties}. Further,  a solution of NCP with neural networks algorithm was presented in  \cite{liao2001solving}. 
Among
these algorithms, the \textit{modified Newton algorithm} is preferred.

The modified Newton algorithm is a reformulation of
 system  of differential equations  \eqref{3} that converts it to a non-modal
optimisation problem which is convergent under some conditions, see Section~8.7~of~\cite{bazaraa2013nonlinear}. This reformulation is given by the system of equations:
\begin{equation}
\label{4}
%\phi(x_i,F(x_i))=0, \qquad\quad \forall i=1, 2, ..., n\,,\\
%\vspace{2cm}
%&
\phi(x)=\left(\begin{array}{lcl}\phi_1(x)\\ \vdots\\
\phi_n(x)\end{array}\right)=
\left(\begin{array}{lcl}\phi(x_1,F_1(x))\\ \vdots\\
\phi(x_n,F_n(x))\end{array}\right)=0,
\end{equation}
where the NCP-function $\phi$  is the  \textit{Fisher-Burmeister NCP-function} (\cite{fischer1992special, fischer1997solution}) defines as follows, 
\begin{equation}\label{5}
\begin{array}{ll}
\phi:I\!\!R^2 \rightarrow I\!\!R, \qquad
\phi(a,b)=\sqrt{a^2+b^2}-(a+b),
\end{array}
\end{equation}
where $\phi(a,b)=0 \Leftrightarrow a\geq 0,\, b\geq0,\, ab=0$.
The resulting system of nonlinear equations (\ref{3}) is
semi-smooth and the NCP has an solution $x$ if and only if $x$ is a solution
for \eqref{4}. Note that the function $\phi$ is locally Lipschitz
and continuous at every point. Thus, the Clarke generalised
Jacobian at every point is defined. In \cite{geiger1996resolution},
%Geiger and Kanzow 
the NCP is recast as an unconstrained minimisation problem and 
the natural merit function $G:I\!\!R^n \rightarrow I\!\!R$ given by
\begin{equation}\label{7} G(x)=\frac{1}{2}
\sum_{i=1}^{n}\phi(x_i,F_i(x))^2=\frac{1}{2}\|\phi(x)\|^2,
\end{equation}
is considered. Also, it is proved that any stationary point of the above function is a
solution of the NCP if the mapping $F$ involved in  Eq.~\eqref{3} is
continuously differentiable and monotone. Hence, the solution of
the NCP are given by solving
\begin{equation}\label{8}
\min_{x\in I\!\!R^n} G(x).
 \end{equation}
 The NCP mainly applied in solving the \textit{traffic equilibrium problem}
which is studied here, and the \textit{Karush-Kuhn-Tucker (KKT) optimality conditions}
 \cite{gordon2012karush}.
%%%%%%%%%%%%%%%%%%
%%%%%%%%%%%%%%%%%%
%%%%%%%%%%%%%%%%%%
\subsection{Traffic equilibrium problem}
%Equilibrium usually interprets as  a state of equality between demand and supply.
Network equilibrium
models  have a variety of applications in urban transportation,
energy distribution, game theory, spatially separated economic
markets, electrical networks, and water resource planning, \cite{aashtiani1981equilibria}.  Although the idea of traffic equilibrium originated as early as
1924 with Frank Knight \cite{knight1924some}, the mathematical description of equilibrium was provided for the first time in 1952 by Wardrop.  
The \textit{Wardrop's equilibrium principals} \cite{wardrop1952road} are  related to the idea of \textit{Nash
equilibrium} \cite{li2009existence} in game theory which is  developed separately.
However, analysing the transportation networks with many players involved is more difficult than analysing games with a few  number  of
players. The Wardrop's principals states as follows:
\paragraph{Wardrop's first principle}
The travel times in all routes actually used are equal and less than those which would be experienced by a single vehicle on any unused route. Each user
non-cooperatively seeks to minimise their transportation cost/time. A traffic flow that satisfies this principle is usually named
 \textit{user equilibrium (UE) flow}.
% Specifically, a user-optimised equilibrium is reached when no user may lower his transportation cost through unilateral action.

\paragraph{Wardrop's second principle} At equilibrium the
average journey time is minimum. This implies that each user behaves
cooperatively in choosing their route to ensure the minimum cost of the whole system. A traffic flow satisfying
Wardrop's second principle is  generally referred to as \textit{system optimal (SO) flow}.

%Wardrop's first principle of route choice, which is identical to the
%notion postulated by Knight, became accepted as a sound and simple
%behavioural principle to describe the spreading of trips over
%alternate routes due to congested conditions. 
 %
Wardrop's principles can be presented mathematically as follows. 
Assume that $R_{ab}$ denotes the set of simple (loop-free) routes for
the origin-destination pair $(a,b)$, $h_r$ denotes the flow on the route $r\in R_{ab}$, and $c_r$ denotes the travel time on the route $r \in R_{ab}$ as experienced by
an individual user. If the flow does  not
depend on time, i.e. in a static model (against the dynamic
model), the Wardrop's principles are written as
\begin{equation}\label{9}
 h_r^T(c_r(h_r)-\pi_{ab})=0,\quad
   c_r(h_r)-\pi_{ab} \geq 0,\quad h_r=d_r(\pi_{ab}),
 \quad h_r \geq 0, \quad \pi_{ab} \geq 0,
\end{equation}
where $\pi_{ab}$ is the minimal rout time and $d_r$ is the demand function
of the route $r\in R_{ab}$.

In  \cite{aashtiani1981equilibria},  an equilibrium model for
predicting traffic flow on a congested transportation network is proposed  based on 
using the Wardrop's first principle
%Let $F:I\!\!R^n \rightarrow I\!\!R^n$
%be a continuous function. For then we could study the  NCP by finding $x \geq
%0$, such that $F(x)^T.x=0$ and $F(x)\geq 0$. 
and write it  as an NCP (Eq.~\eqref{3}) in the following form:
\begin{equation}\label{10}
\displaystyle
\begin{array}{cc}
F(x)^T \cdot x=\Big(c_r(h_r)-\pi_{ab} \quad h_r-d_r(\pi_{ab})\Big)\cdot \left(\begin{array}{cc}
h_r\\
\pi_{ab}
\end{array}\right)=0, \\
x=\left(\begin{array}{cc}
h_r\\
\pi_{ab}
\end{array}\right)
\geq 0,\qquad  
F(x)=\left(\begin{array}{cc}
c_r(h_r)-\pi_{ab}\\
h_r-d_r(\pi_{ab})\end{array}\right)\geq0,
\end{array}
\end{equation}
where $F:I\!\!R^n \rightarrow I\!\!R^n$ is a continuous function.
For reducing this NCP to an
unconstrained global minimisation problem,  
gap function $G$ which is smooth, convex and bounded is applied,  see   \cite{facchinei1997new}:
\begin{equation}\label{11}
 G(x)=\sum^n_{i=1}
\varphi \big(x_i,F(x_i)\big),
\end{equation}
 where   $\varphi=\frac{1}{2}\phi^2$ and  $\phi$ is the
Fischer-Burmeister function given by (\ref{5}). So,  solving
the NCP \eqref{10}   is equivalent to finding a general unconstrained solution for Eq.~(\ref{8}) with $G$ given in \eqref{11}.

%%%%%%%%%%%%%%%%%%%%%%
%%%%%%%%%%%%%%%%%%%%%%
%%%%%%%%%%%%%%%%%%%%%%
\section{A Finsler geometric model of traffic problem}
\label{Sec:Finsler-Traffic-model}
 In this section, we apply geometric methods to provide a mathematical model for the transportation network. Unless stated otherwise we will assume further that 
 \begin{enumerate} 
 \item[I] All origin-destination pairs are distinct.
 \item[II] The network is strongly connected, that is, at least one route joins each origin-
destination pair.
\item [III] The transport demand is either  a constant number or a function of  travel cost/time.
\item [IV] All travellers have perfect information about their
travel. Therefore, we have a deterministic model, unlike a
stochastic model where travellers choose their paths randomly.
\end{enumerate}

For any given pair of the origin-destination points in the vehicle network, we usually have different route options. 
 These routes are made by a string of streets and cross-sections and for each route we have an estimated travel time.
This time depends on some factors such as the capacity of
streets, the number of stops behind cross-sections and the travel demand of that specific route.
 
% The traffic congestion is a road condition characterised by slower speeds, longer trip times and increased queuing. It occurs when
%roadway (or path) demand is greater than its capacity.

Now assume that we have  a  Riemannian metric $g$ on  $I\!\!R^2$ and
 its  corresponding norm-squared  for tangent vectors
$y\in T_x I\!\!R^2$ is  given by $\|y\|_g^2:=g(y,y)$. In traffic modelling, we can interpret it as the time it takes, using
a car with fixed power, to travel from the base point of the
vector $y$ to its tip.

Let $u\in T_x I\!\!R^2$ be a unit vector and denote the
traffic congestion or any external factor that increases the traffic
congestion by a vector $\omega\in T_x I\!\!R^2$ such that $\|\omega\|_g<1$.
Before $\omega$ sets in, a journey from the base to the tip of any
$u$ would take one unit of time. After the  effect of $\omega$, within the
same one unit of time, we traverse not $u$ but 
$v=u-\omega$ instead.  This is because  traffic congestion is a vector in
the opposite direction of flow. The measure of this new vector is not
equal to 1 ($\|v\|_g\neq1$). This argument for using vector's length is the key idea of a technique  known as \textit{Okubo's technique} \cite{antonelli2013theory}. In fact, the geometry of movement in the former case
(without considering the external factor) is the  Riemannian
geometry, rather than the Finslerian geometry in the latter case.
 Here, we consider the effect of traffic
congestion (the vector field $\omega$) and introduce a new metric
$\overline{F}$ such that  $v$ be a unit vector with respect to the new norm corresponding to $\overline{F}$.
\begin{theorem}\label{Randers}
Let $g$ be  a Riemannian metric and $\omega$ a vector field
on $I\!\!R^2$ such that $\|\omega\|_g<1$. If $\omega$ indicates the
traffic congestion, then the travel time  for a car with a fixed
power to pass through a vector field $y\in T_x I\!\!R^2$  is measured by
the Randers metric $\overline{F}(x,y)=\frac{\|y\|_g}{1-\|\omega\|_g}$. 
\end{theorem}
\begin{proof}
Let $u\in T_x I\!\!R^2$ be a unit vector. After taking the effect of the vector
$\omega$ into account,  we have $v=u-\omega$, where $\|v\|_g=1-\|\omega\|_g$. Now, a  Finsler structure $\overline{F}$ that satisfies $\overline{F}(v)=1$ would be 
 $\displaystyle \overline{F}(v)=\frac{\|v\|_g}{1-\|\omega\|_g}$. For an arbitrary vector field $y\in T_x I\!\!R^2$, the corresponding Finsler metric $\overline{F}$ is 
$$\overline{F}(x,y)=\frac{\|y\|_g}{1-\|\omega\|_g}.$$
 If we put $\lambda=1-\|\omega\|_g^2$ then $\overline{F}(x,y)$ can be written as
\begin{equation}\label{eq:Randers}
\overline{F}(x,y)=\frac{\|y\|_g}{\lambda}+\frac{\|\omega\|_g\|y\|_g}{\lambda}=\frac{\sqrt{g_{ij}
y^i y^j}}{\lambda}+ \frac{g_{ij}\omega^i y^j}{\lambda},
\end{equation}
which is obviously a Randers metric.
%, i.e. $(\alpha+\beta)$-metric, with $\alpha=\frac{\sqrt{g_{ij} y^i y^j}}{\lambda}$ and $\beta=b_jy^j=\frac{\omega_j}{\lambda}y^j$, where$\omega_j=g_{ij}\omega^i$.  
 \end{proof}
 Therefore,  the optimal/shortest route between any
 arbitrary origin-destination pair $(p,q)$ is a  geodesic of the corresponding  Randers metric passing through these points . Equivalently, we need to find
minimums of the following integral
\begin{equation}\label{2}
J_{\overline{F}}(C)=\int_0^1\overline{F}\big(x(t),y(t)\big)dt, \qquad y(t)=\frac{dx}{dt},
\end{equation}
where $\overline{F}$ is given by~\eqref{eq:Randers} and $C:[0,1]\rightarrow I\!\!R$ defined by $t\rightarrow x(t)$  such that $x(0)=p$ and $x(1)=q$.  It can be shown that by assuming
$L=\frac{\overline{F}^2}{2}$, any  minimal point $x(t)$ of the length
integral~\eqref{2} is a solution of the following system
of differential equations, see Chapter 3 of \cite{shen2001lectures}:
\begin{equation}\label{eq:eq}
L_{ij}\frac{d^2x^i}{dt^2}+\frac{\partial L_i}{\partial x^j}\frac{dx^j}{dt}-\frac{\partial L}{\partial x^i}=0,
\end{equation}
where $L_i=\frac{\partial L}{\partial y^i}$,
$L_{ij}=\frac{\partial^2 L}{\partial y^i\partial y^j}$. Therefore, 
\begin{corollary}
If we consider the traffic congestion as an external factor  that effects on vehicle transition network, then the differential equation of the shortest route  between any origin-destination pair is given by Eq.~\eqref{eq:eq}.
\end{corollary}

 %%%%%%%%%%%%%%%%%%%%%%%
 %%%%%%%%%%%%%%%%%%%%%%%
 %%%%%%%%%%%%%%%%%%%%%%%%
 \section{A Finslerian model for Wardrop's user equilibrium problem}
 \label{Sec:Ward}
 In this section we show that  the Fischer-Burmeister function is a
special Randers metric, and therefore, the equations of its geodesics are equations of the minimal paths. Then, we present a Finsler geometrical model for solving the dynamical Wardrop's user equilibrium problem that may applied in solving different equilibrium problems.
% In section 3, we saw
%that if we consider the traffic congestion as an external factor that
%effects  traffic flow, then the geometrical method shows that the
%geodesics of a Randers metric are the shortest path for travelling from
%an origin point to a destination point in a transition network. In
%section 4, we proved that the Fischer-Burmeister function is a
%special Randers metric, and therefore, the equations of the minimal path
%are its geodesics. Now regarding what presented in the previous
%section, we give a mathematical model for solving the dynamical
%Wardrop's user equilibrium problem, by using Finsler geometry.
 \subsection{Finsler geometry and traffic equilibrium}
 \label{Sec:Finsler-Traffic-equilib}
 Here, we illustrate that the Fischer-Burmeister function is a special
Randers metric.
 \begin{proposition} \label{Fi_Bu}
  The Fischer-Burmeister function
is a Finsler metric of Randers type  and its geodesics are
the optimised paths in the  NCP.
 \end{proposition}
\begin{proof}
Let $g$ be a Euclidean metric on $I\!\!R^2$ given by
$g(a,b)=\sqrt{a^1b^1+a^2b^2}$, where $a=(a^1,a^2)$ and $b=(b^1,b^2)$
are two arbitrary points on $I\!\!R^2$. The local coordinate system $(x^1,x^2)$ on $I\!\!R^2$  induces the local
 fields of frames $\left(\frac{\partial}{\partial x^1},
\frac{\partial}{\partial x^2}\right)$ and $(dx^1,dx^2)$ on the tangent
space $T_x I\!\!R^2$ and its dual $T^*_x I\!\!R^2$.
Therefore, we have $g\left(\frac{\partial}{\partial x^1},
\frac{\partial}{\partial x^2}\right)=\sqrt{(dx^1)^2+(dx^2)^2}$ and from Eq.
(\ref{5}), the Fischer-Burmeister function on $I\!\!R^2$  is given by:
\begin{equation}\label{eq:FB}
 \phi(dx^1,dx^2)=\sqrt{(dx^1)^2+(dx^2)^2}-(dx^1+dx^2 ).
 \end{equation}
 The first term in the right hand side of \eqref{eq:FB} can be written as $\sqrt{\delta_{ij}dx^idx^j}$, where $\delta_{ij}$ is the
Kronecker symbol and $i, j$ run over the range 1,2. So,  it  is a Euclidean
(or Riemannian) metric on $I\!\!R^2$.
 The second term in \eqref{eq:FB} is a
differential 1-form $\beta=\omega_idx^i$, where $\forall i: \omega_i=-1$. Therefore,
the Fischer-Burmeister function is a Randers metric.
 \end{proof}
 By virtue of the above proposition, solutions of the NCP, i.e. the optimised paths of the Fischer-Burmeister function, are  geodesics of a  Randers
metric given by Eq.~\eqref{eq:eq}.
\subsection{The dynamical transition network}
\label{Sec:dynamic}
%%%%%%%%%%%%%%%%%%%
The traffic congestion and transport demand in some hours of a day
(peak hours) are   extremely increasing. A typical solution to address this issue,  instead of building additional infrastructure, is introduction of dynamic elements to the road traffic management. This includes for instance using sensors  to measure traffic flows and automatic interconnected guidance systems
(for example traffic signs which open a lane in different directions
depending on the time of day) to manage traffic in peak
hours. Here, we generalise the method of solving the NCP (in
the previous subsection) to dynamic systems.
\begin{theorem}\label{dynamic}
If the traffic flow  and the travel time  are functions of time $t$,
then the solutions of traffic equilibrium problem are minimal of the
following integration
\begin{equation}\label{12}
G(t)=\int^{t_1}_{t_0}\varphi (h(t),c(t)) h'(t) dt,
\end{equation}
 where $\varphi=\frac{1}{2}\phi$ and $\phi$ is the Fischer-Burmeister  function, $h(t)$ is  the traffic flow,
   $c(t)$ is the travel time, $t_0$ is the start of travel time, and $t_1$ is the end of travel time.
\end{theorem}
\begin{proof}
Let the traffic flow  $h(t)$ and the
travel time $c(t)$ (which is in general a function of $h$) be
functions of time and assume that the minimal rout time  $\pi$ in \eqref{9} is zero. Then, using the Fischer-Burmeister  function $\phi$  on $I\!\!R^2$ (Eq. \eqref{eq:FB}),  the Eq. (\ref{7}) can be
written as
$$G(h)=\int^{h_1}_{h_0}\varphi (h,c(h))  dh,$$
 where $\varphi=\frac{1}{2}\phi$. Since
  $dh=h'(t)dt$, this
equation is equivalent to Eq. (\ref{12}). Now  by virtue of a well-known fact in the calculus of variation the minimal of integration
(\ref{12}) are solutions of traffic equilibrium problem or present the shortest route (the route which needs the
shortest time duration to passing through it).
 \end{proof}
 %%%%%%%%%%%%%%%%%%%%%
 %%%%%%%%%%%%%%%%%%%%%%
 \subsection{The Finsler geometrical model}
 Now,  we can present a  Finsler geometrical model for solving the dynamical Wardrop's user equilibrium
problem: 
% In section 3, we saw
%that if we consider the traffic congestion as an external factor that
%effects  traffic flow, then the geometrical method shows that the
%geodesics of a Randers metric are the shortest path for travelling from
%an origin point to a destination point in a transition network. In
%section 4, we proved that the Fischer-Burmeister function is a
%special Randers metric, and therefore, the equations of the minimal path
%are its geodesics. Now regarding what presented in the previous
%section, we give a mathematical model for solving the dynamical
%Wardrop's user equilibrium problem, by using Finsler geometry.
 \begin{lemma} \label{Finsler}
Any  Finsler  structure $\overline{F}(x,y)$  on $I\!\!R^2$ 
  can be considered as a gap function in solving the  NCP.
 \end{lemma}
\begin{proof}
Assume that the transit network has all conditions I-IV in Section~\ref{Sec:Finsler-Traffic-model}  and also traffic
flow and the other ingredients in the transit network are functions of
time. Further assume that $\overline{F}(x,y):T I\!\!R^2 \backslash \{0\}\rightarrow I\!\!R^+$
is a Finsler metric on $I\!\!R^2$ and a car with a fixed
power travels from point $a$ to point $b$ trough a smooth curve $r$ on $I\!\!R^2$.  From the calculus of variation the minimal of the following
integration provide the shortest curve between points $a$ and $b$:
\begin{equation}\label{13}
L_{\overline{F}}(r)=\int_a^b\overline{F}(x(t),y(t))dt,
\end{equation}
where $y(t)=\frac{dx}{dt}$. \\
Since $\overline{F}$ is a Finsler
    structure, the following conditions are satisfied (see for instance\cite{antonelli2013theory, shen2001lectures}):

 1. The length of the curve $r$ is independent of the parameter $t$.

 2. The length of the curve $r$  is always a positive number.

 3. If we put
   $(g_{ij}):=\left({1 \over 2} \left[ \frac{\partial^{2}}{\partial y^{i}\partial y^{j}}\overline{F}^2 \right]\right)$,
    then  we have $\det(g_{ij})\neq0$.

    4. $\overline{F}(x,y)$ is a differentiable function.\\
Therefore, $\overline{F}(x,y)$ meets all the conditions of a gap function
in solving the NCP, see \cite{bazaraa2013nonlinear}.
 \end{proof}
\begin{theorem} \label{Optimise}
The optimised routes for the traffic equilibrium problem in a
dynamical transition network are geodesics of a Finsler structure
$\overline{F}(x,y)$ of Randers type, where $x(t)$ denotes the
traffic flow at time $t$ and $y(t)=\dot{x}(t)$ is the velocity of a
car with fixed power at time $t$.
\end{theorem}
\begin{proof}
Let $x(t)$ be the traffic flow and $y(t)=\dot{x}(t)$ the velocity of
a car with fixed power at time $t$. Then, the Wardrop's user
equilibrium principle can be written as
$$\overline{F}(x,y)=0\Leftrightarrow x\geq0,\quad y\geq0,\quad x.y=0.$$
Here, $\overline{F}(x,y)=0$ implies that the traffic congestion is
extremely increased and the car velocity is zero ($y=0$) or the car stops, so we have $x.y=0$. In Section~\ref{Sec:Finsler-Traffic-model}, we
show that the geometrical method of solving the traffic problem
leads to  finding the geodesics of a  Randers metric. On the other hand, Proposition~\ref{Fi_Bu} implies that the common applied mathematical method, i.e. finding the minimums of the Fischer-Burmeister function,  can be seen as a special case of this geometrical method. Thus, if the given Finsler metric in Lemma \ref{Finsler}be a Randers metric, it is a solution  of the NCP.
Therefore, from Theorem \ref{dynamic}, the
equations of optimised routes for the traffic equilibrium problem in
a dynamical transition network are given by the geodesic equations
of a Finsler metric of Randers type.
 \end{proof}

\bibliography{biblio}
\end{document}